\documentclass[11pt,leqno]{amsart}
\usepackage[utf8]{inputenc}
\DeclareTextSymbolDefault{\textsection}{OMS}
\usepackage[a4paper,left=29.7377mm,textwidth=145mm,
  top=46.3844mm,textheight=210mm,headheight=8pt,
  headsep=10mm,footskip=12pt]{geometry}
\usepackage{etoolbox}
\makeatletter
\patchcmd{\@setabstracta}{\skip@20\p@}{\skip@44.13\p@}{}{}
\makeatother
\usepackage{amsmath,amssymb,mathtools}
\usepackage{enumitem}
\usepackage{needspace}
\usepackage[hidelinks,pdfencoding=auto]{hyperref}
\usepackage[capitalise,nameinlink,noabbrev]{cleveref}
\hypersetup{pdftitle={THE FIRST POSITIVE EIGENVALUE OF THE DIRAC OPERATOR ON ISOPARAMETRIC HYPERSURFACES IN THE UNIT SPHERE},pdfauthor={Jianquan Ge and Xize Yang}}
\newtheorem{theorem}{Theorem}[section]
\newtheorem{proposition}[theorem]{Proposition}
\newtheorem{lemma}[theorem]{Lemma}
\newtheorem{corollary}[theorem]{Corollary}
\newtheorem{conjecture}[theorem]{Conjecture}
\theoremstyle{definition}

\theoremstyle{remark}

\numberwithin{equation}{section}
\newcommand{\R}{\mathbb R}
\newcommand{\C}{\mathbb C}
\newcommand{\Z}{\mathbb Z}
\newcommand{\Id}{\mathrm{Id}}
\newcommand{\Scal}{\operatorname{Scal}}
\newcommand{\Spec}{\operatorname{Spec}}
\newcommand{\tr}{\operatorname{tr}}
\newcommand{\rank}{\operatorname{rank}}
\newcommand{\vol}{\,dv}
\newcommand{\DL}{D_L}
\newcommand{\dd}{\,d}
\newcommand{\ip}[2]{\langle #1,#2\rangle}
\newcommand{\HS}{\mathcal H}
\DeclareMathOperator{\Spin}{Spin}
\DeclareMathOperator{\SO}{SO}

\newcommand{\E}{\widehat{\mathcal E}}
\newcommand{\De}{\widehat{\mathcal D}}
\newcommand{\LA}{\mathcal L_A}
\newcommand{\KK}{\mathsf K}
\newcommand{\NN}{\mathsf B}
\newcommand{\Pe}{P_{\mathrm e}}
\newcommand{\volM}{\operatorname{vol}_M}
\newcommand{\norm}[1]{\lVert #1\rVert}
\newcommand{\Ltwo}[2]{(#1,#2)}
\DeclareMathOperator{\Rea}{Re}

\allowdisplaybreaks[2]
\title[First positive Dirac eigenvalue]{THE FIRST POSITIVE EIGENVALUE OF THE DIRAC OPERATOR\\
ON ISOPARAMETRIC HYPERSURFACES IN THE UNIT SPHERE}
\author[J. Q. Ge]{Jianquan Ge}
\address{School of Mathematical Sciences, Beijing Normal University,
Beijing 100875, P. R. China}
\email{jqge@bnu.edu.cn}
\author[X. Yang]{Xize Yang$^{*}$}
\address{Chern Institute of Mathematics and LPMC, Nankai University,
Tianjin 300071, P. R. China}
\email{xize.yang@nankai.edu.cn}
\thanks{$^{*}$Corresponding author.}
\thanks{J. Ge is partially supported by the NSFC (No. 12571049)
and the Fundamental Research Funds for the Central Universities.}
\date{}

\begin{document}
\begin{abstract}
We prove that the intrinsic Dirac spectrum of any
minimal isoparametric hypersurface $M^n\subset S^{n+1}(1)$, equipped with
the induced spin structure, is disjoint from $(-n/2,n/2)$ and contains
both endpoints. In particular, its first positive Dirac eigenvalue is $n/2$ and the proof is independent of the classification of isoparametric hypersurfaces.
We also show that, for any closed connected embedded minimal hypersurface
in the unit sphere with its induced spin structure, this Dirac spectral gap
implies that the first positive
Laplace eigenvalue is $n$. These two results give an alternative spinorial
proof of Tang and Yan's theorem and show that the same Dirac spectral gap would
suffice for the general Yau conjecture.
\end{abstract}
\maketitle

\section{Introduction}

Yau's conjecture \cite[Problem~100]{YauProblem} asserts that every closed
connected minimal hypersurface $M^n$ embedded in the unit sphere
$S^{n+1}(1)$ satisfies
\[
 \lambda_1(\Delta_M)=n,
\]
where $\Delta_M=-\operatorname{div}\nabla$ is the nonnegative
Laplace--Beltrami operator. By Takahashi's formula \cite{Takahashi}, the
coordinate functions of the inclusion satisfy $\Delta_M x_j=nx_j$.
Thus $n$ is an eigenvalue, and the problem is to prove the lower bound
$\lambda_1(\Delta_M)\ge n$. The classical estimate of Choi and Wang
\cite{ChoiWang} gives $\lambda_1(\Delta_M)\ge n/2$.

Tang and Yan \cite[Theorem~1.1]{TangYan} proved the conjecture for all
minimal isoparametric hypersurfaces in 2013. Their theorem built on
earlier results for homogeneous and inhomogeneous families due to Muto,
Ohnita, and Urakawa \cite{MutoOhnitaUrakawa1984}, Kotani
\cite{Kotani1985}, Muto \cite{Muto1988}, and Solomon
\cite{Solomon1990I,Solomon1990II}. Tang, Xie, and Yan
\cite[Theorem~1.1]{TangXieYan2014} subsequently gave a direct proof for
six principal curvatures of multiplicity two, without using the
classification of that family. We refer to
\cite[\S4.3]{GeQianTangYanSurvey} for an account of these developments.

In this paper, we study Yau's conjecture through the Dirac operator.
Let $D_M$ denote the intrinsic Dirac operator for the spin structure
induced from the sphere, and let $\lambda_1^+(D_M)$ be its first positive
eigenvalue. On a closed embedded minimal hypersurface, ambient Killing spinors give the
eigenvalues $\pm n/2$; a construction valid in every dimension is recalled
in \cref{sec:killing}. This suggests the spinorial question of whether
the interval $(-n/2,n/2)$ contains any Dirac eigenvalues. Our main results
establish this spectral gap for minimal isoparametric hypersurfaces and
show that, for an arbitrary embedded minimal hypersurface, the gap implies
$\lambda_1(\Delta_M)=n$.

We recall several intrinsic and extrinsic estimates for Dirac eigenvalues. For a closed connected embedded hypersurface in the unit sphere, let
$A$ be the shape operator for a chosen unit normal and let $H=\tr A/n$
be the normalized mean curvature. For the induced spin structure,
B\"ar's upper bound \cite{BaerExtrinsic} reads
\begin{equation}\label{eq:baer}
 \min\Spec(D_M^2)
 \leq \frac{n^2}{4\operatorname{Vol}(M)}
       \int_M(1+H^2)\vol_M.
\end{equation}
Equality holds for round hyperspheres and, by work of Ginoux, for products
of round spheres \cite{GinouxProducts} and the three-dimensional minimal
Cartan hypersurface \cite{GinouxCartan}. Ginoux anticipated equality for
isoparametric hypersurfaces in \cite[\S1.2.3, p.~13]{GinouxHabil}.
Ginoux, Habib, and Raulot \cite[p.~81]{GinouxHabibRaulot} also described
the spherical equality problem as a possible spinorial analogue of Yau's
conjecture.

Classical intrinsic estimates follow from the Schr\"odinger--Lichnerowicz formula \eqref{eq:SL} \cite{Lichnerowicz}.
The first is Friedrich's intrinsic lower bound \cite{Friedrich}
\[
 \min\Spec(D_M^2)\ge \frac{n}{4(n-1)}\inf_M\Scal_M,
 \qquad n\ge2.
\]
Hijazi obtained refinements using conformal geometry for $n\ge3$
\cite{HijaziConformal} and the energy--momentum tensor
\cite{HijaziEnergyMomentum}. For a minimal hypersurface in the unit
sphere, however, the Gauss equation gives
$\Scal_M=n(n-1)-|A|^2$. The term $-|A|^2$ therefore prevents
Friedrich's estimate alone from reaching $n^2/4$ unless the hypersurface
is totally geodesic.

Extrinsic lower bounds provide a different source of estimates.
Zhang's inequalities \cite{Zhang1998,Zhang1999} and their conformal
refinements by Hijazi and Zhang \cite{HijaziZhangI} involve the scalar
curvature, mean curvature, and energy--momentum tensor for operators
defined through the ambient spin connection. For the intrinsic boundary
Dirac operator, the spinorial Reilly estimate of Hijazi, Montiel, and
Zhang \cite{HMZEmbedded} gives a lower bound in terms of the boundary
mean curvature. That bound vanishes on a minimal boundary; their
conformal estimate \cite{HMZConformal} also involves the geometry of the
filling.\footnote{Related spherical bounds for spin fillings are discussed
in Montiel's preprint \cite{Montiel2024}. None of those bounds is used
in the proofs below.}
These estimates motivate a lower bound that uses the geometry of the
principal distributions directly.

Our first theorem gives the sharp Dirac spectral gap in the isoparametric case.
Throughout the paper, \emph{closed} means compact without boundary, hypersurfaces are closed, connected and embedded with the induced spin structure from the standard one on the ambient sphere.

\Needspace{10\baselineskip}
\begin{theorem}\label{thm:minimal}
Let $n\geq1$, and let $M^n\subset S^{n+1}(1)$ be a closed minimal isoparametric hypersurface. For the induced spin structure,
its intrinsic Dirac operator satisfies
\begin{equation}\label{eq:main}
 \Spec(D_M)\cap(-n/2,n/2)=\varnothing,
 \qquad
 \{-n/2,n/2\}\subset\Spec(D_M).
\end{equation}
Consequently,
\[
 \lambda_1^+(D_M)=\frac n2,
 \qquad
 \min\Spec(D_M^2)=\frac{n^2}{4}.
\]
\end{theorem}

The proof rests on the following estimate for the spin connection,
which holds without a minimality assumption.

\begin{theorem}\label{thm:connection}
Let $n\geq1$, and let $M^n\subset S^{n+1}(1)$ be a closed isoparametric hypersurface with the induced spin structure.
Every smooth spinor $\psi$ on $M$ satisfies
\begin{equation}\label{eq:connection-main}
 \int_M|\nabla^\Sigma\psi|^2\vol_M
 \geq \frac{n+|A|^2}{4}\int_M|\psi|^2\vol_M.
\end{equation}
The constant $(n+|A|^2)/4$ is optimal.
\end{theorem}

For a minimal hypersurface, the term $|A|^2$ in
\eqref{eq:connection-main} cancels the corresponding term in the Gauss
equation when the estimate is inserted into the integrated
Schr\"odinger--Lichnerowicz formula \eqref{eq:integrated-SL}. This gives the lower bound in
\cref{thm:minimal}; ambient Killing spinors realize both endpoints.
The same argument applies to nonminimal isoparametric hypersurfaces
and proves equality in B\"ar's estimate \eqref{eq:baer}, confirming Ginoux's conjecture.

\begin{corollary}\label[corollary]{cor:cmc}
Let $M^n\subset S^{n+1}(1)$ be as in \cref{thm:connection}, and let
$H=\tr A/n$. 
Then
\[
 \Spec(D_M)\cap(-\frac n2\sqrt{1+H^2},\frac n2\sqrt{1+H^2})=\varnothing,\quad
 \{\pm\frac n2\sqrt{1+H^2}\}\subset\Spec(D_M).
\]
In particular,
\[
 \lambda_1^+(D_M)=\frac n2\sqrt{1+H^2},
 \qquad
 \min\Spec(D_M^2)=\frac{n^2}{4}(1+H^2).
\]
\end{corollary}

Our next theorem explains how the Dirac spectral gap leads to Yau's
conjecture. Its hypothesis is spectral, and no
isoparametric assumption is needed.

\begin{theorem}\label{thm:laplace-transfer}
Let $M^n\subset S^{n+1}(1)$, $n\ge1$, be a closed connected
embedded minimal hypersurface with the induced spin structure. If $\min\Spec(D_M^2)=\frac{n^2}{4}$, or equivalently,
\begin{equation}\label{eq:transfer-gap}
 \Spec(D_M)\cap(-n/2,n/2)=\varnothing,
\end{equation}
then $\lambda_1(\Delta_M)=n$.
\end{theorem}

Combining \cref{thm:minimal,thm:laplace-transfer} recovers the theorem
of Tang and Yan.

\begin{corollary}[Tang--Yan \cite{TangYan}]\label[corollary]{cor:tang-yan}
Every closed minimal isoparametric hypersurface
$M^n\subset S^{n+1}(1)$ satisfies $\lambda_1(\Delta_M)=n$.
\end{corollary}

The Dirac spectral gap in \cref{thm:minimal} is proved uniformly for all isoparametric hypersurfaces, using only their common geometric properties without referring to the homogeneity or the classification. Consequently, \cref{thm:minimal,thm:laplace-transfer} give an
alternative spinorial proof of the isoparametric case of Yau's conjecture.

The proof uses the full interval exclusion in \eqref{eq:transfer-gap}.
The intrinsic Dirac spectrum need not be symmetric
about zero \cite[Theorem~1.3.7 and Corollary~2.1.5]{GinouxSpectrum},
so a bound on the first positive eigenvalue alone does not in general
control the negative spectrum. Our positive lower bounds for $D_M^2$ control both sides of the spectrum
and exclude zero at the same time, without a distinction between dimensions. This motivates the following spinorial version of Yau's conjecture.

\begin{conjecture}\label{conj-spinorial-Yau}
 Let $M^n\subset S^{n+1}(1)$ be a closed connected
embedded minimal hypersurface with the induced spin structure. Then the
smallest eigenvalue of the square of its intrinsic Dirac operator is $\frac{n^2}{4}$.
\end{conjecture}

By \cref{thm:laplace-transfer}, Conjecture \ref{conj-spinorial-Yau} implies the general Yau
conjecture and it is proved in the isoparametric case in \cref{thm:minimal}. 

We conclude the introduction with the structure of the proof. The
geometric input is the curvature-sphere geometry developed by Cartan
\cite{Cartan1938} and M\"unzner \cite{MunznerI,MunznerII}. It applies
equally to the inhomogeneous examples of Ozeki and Takeuchi
\cite{OzekiTakeuchi1975} and those obtained from Clifford systems by
Ferus, Karcher, and M\"unzner \cite{FKM1981}; see
\cite[\S\S1--2]{GeQianTangYanSurvey},
\cite[Chapters~2--3]{CecilRyan2015}, and \cite{CecilRyan2025} for
background. Sections~\ref{sec:leaf-estimates}--\ref{sec:global} establish the
spin connection estimate on each curvature sphere and integrate it over
$M$ using the focal maps. Section~\ref{sec:killing} constructs the
endpoint eigenspinors and completes the Dirac estimates.
The passage to the Laplacian uses an ambient Clifford description of the
Dirac operator; see \cite{Trautman} for related hypersurface constructions.
In Section~\ref{sec:laplace}, a Clifford trace identity converts the
Dirac gap into an estimate for a first-order operator on
differential forms. After identifying its kernel, we use the Fredholm
alternative to obtain the lower bound for the scalar Laplacian.

\section{Preliminaries}\label{sec:conventions}

Our conventions follow \cite[\S\S1--2]{GeQianTangYanSurvey} and
\cite{Cartan1938,CecilRyan2015,CecilRyan2025,MunznerI,MunznerII} for isoparametric geometry, and
\cite[Chapter~II, \S\S1--5]{LawsonMichelsohn} and
\cite[\S\S2.1, 2.3--2.5]{BourguignonEtAl} for spin geometry.

Let $x:M^n\hookrightarrow S^{n+1}(1)\subset\R^{n+2}$, $n\geq1$,
be a closed connected embedded hypersurface. Embeddedness assures that $M$ admits a smooth spherical unit normal $\nu$. Fix such a normal
and use it to orient $M$. Write $\nabla^{\R}$, $\bar\nabla$, and
$\nabla$ for the Euclidean, spherical, and intrinsic Levi--Civita
connections, respectively. Our shape operator convention is
\begin{equation}\label{eq:shape}
 A(X)=-\bar\nabla_X\nu=-\nabla_X^{\R}\nu,
 \qquad H=\frac1n\tr A.
\end{equation}
Thus $H$ is the normalized mean curvature, and the Gauss formulas are
\begin{align}
 \bar\nabla_XY&=\nabla_XY+\ip{AX}{Y}\nu,
       \label{eq:spherical-gauss}\\
 \nabla_X^{\R}Y
 &=\nabla_XY+\ip{AX}{Y}\nu-\ip{X}{Y}x.
       \label{eq:euclidean-gauss}
\end{align}
With the convention
$R(X,Y)=\nabla_X\nabla_Y-\nabla_Y\nabla_X-\nabla_{[X,Y]}$,
the Gauss and Codazzi equations read
\begin{align}
 R^M(X,Y)Z
 &=\ip{Y}{Z}X-\ip{X}{Z}Y
   +\ip{AY}{Z}AX-\ip{AX}{Z}AY,
       \label{eq:gauss}\\
 (\nabla_XA)Y&=(\nabla_YA)X.
       \label{eq:codazzi}
\end{align}
Tracing the Gauss equation \eqref{eq:gauss}, we express the intrinsic scalar curvature by
\begin{equation}\label{eq:scalar}
 \Scal_M=n(n-1)+(\tr A)^2-|A|^2
       =n(n-1)+n^2H^2-|A|^2.
\end{equation}

A smooth nonconstant function $f$ on the sphere is isoparametric if
\[
 |\bar\nabla f|^2=b(f),\qquad
 \Delta_{S^{n+1}}f=a(f),
\]
for smooth functions $a,b$ of one variable, with
$\Delta_{S^{n+1}}=-\operatorname{div}\bar\nabla$.
The connected components of regular level sets are isoparametric
hypersurfaces. For closed connected embedded hypersurfaces in space forms, this
definition is equivalent to the constancy of the principal curvatures
\cite{Cartan1938,GeQianTangYanSurvey}.

For an isoparametric $M$, write $\kappa_1,\ldots,\kappa_g$ for the
distinct constant principal curvatures, $E_i=\ker(A-\kappa_i\Id)$ for their
principal distributions, and $m_i=\rank E_i$. Each $E_i$ is a smooth
subbundle of constant rank. It is integrable, and its maximal connected
leaves are totally geodesic in $M$. Each such curvature leaf $L$ is a
complete round $m_i$-sphere of radius $r_i=(1+\kappa_i^2)^{-1/2}$, called curvature sphere,
embedded in an affine $(m_i+1)$-plane in $\R^{n+2}$
\cite[Lemma~1.2]{Siffert}. For $X$ tangent to $L$ and $Z$ normal to
$L$ in $M$, the identity $AX=\kappa_i X$ and the Gauss formula
\eqref{eq:euclidean-gauss} give $\nabla_X^{\R}Z=\nabla_X Z$.
It is not hard to prove \cite{CecilRyan2015} that parallel transport in $M$ along $L$
keeps each normal vector constant in $\R^{n+2}$. The normal
spaces of $L$ in $M$ all coincide with a fixed subspace of $\R^{n+2}$,
and any orthonormal basis of this subspace gives a global parallel
normal frame along $L$.

At each point of $M$, the spectral decomposition of $A$ gives
\begin{equation}\label{eq:splitting}
 \begin{gathered}
 TM=E_1\oplus\cdots\oplus E_g\quad\text{orthogonally},
 \qquad \sum_i m_i=n,\\
 |A|^2=\sum_i m_i\kappa_i^2,
 \qquad nH=\sum_i m_i\kappa_i.
 \end{gathered}
\end{equation}
In particular, $H$, $|A|^2$, and $\Scal_M$ are constant.
By M\"unzner's theorem, the principal curvatures can be ordered so that
\[
 \kappa_i=\cot\!\left(\theta+\frac{(i-1)\pi}{g}\right),
 \quad 0<\theta<\frac\pi g,
 \qquad m_{i+2}=m_i,
\]
where the multiplicity indices are cyclic and $g\in\{1,2,3,4,6\}$
\cite{MunznerI,MunznerII}. See also
\cite[Theorem~2.1]{GeQianTangYanSurvey}.

The parallel family and its spherical unit normal are given by
\begin{equation}\label{eq:parallel-family}
 \begin{gathered}
 x_t=\cos t\,x+\sin t\,\nu,
 \qquad \nu_t=-\sin t\,x+\cos t\,\nu,\\
 dx_t=\cos t\,\Id-\sin t\,A.
 \end{gathered}
\end{equation}
Whenever $x_t$ is an immersion, its principal curvatures with respect
to $\nu_t$ are
\[
 \kappa_i(t)=
 \frac{\sin t+\kappa_i\cos t}{\cos t-\kappa_i\sin t}.
\]
For $\theta_i\in(0,\pi)$ with $\cot\theta_i=\kappa_i$, the focal
map $x_{\theta_i}=(\nu+\kappa_i x)/\sqrt{1+\kappa_i^2}$ has
differential with kernel $E_i$. The regular hypersurfaces are tubes
around the two focal submanifolds
\cite{MunznerI,MunznerII,GeQianTangYanSurvey}. We shall estimate the
spin connection leafwise on each curvature sphere and integrate it over
$M$ using the focal maps in Sections~\ref{sec:leaf-estimates}--\ref{sec:global}.

Let $P_{\Spin}(M)$ denote the principal $\Spin(n)$-bundle that
double-covers the oriented orthonormal frame bundle equivariantly
with respect to $\Spin(n)\to\SO(n)$. We equip $M$ with the spin
structure induced from $\R^{n+2}$ through $S^{n+1}$ by the chosen
normals. Equivalently, this structure is induced by the stable
oriented trivialization
\begin{equation}\label{eq:stable-spin}
 TM\oplus\R\nu\oplus\R x\cong M\times\R^{n+2}.
\end{equation}
This choice determines the curvature-circle holonomy in
\cref{lem:circle-holonomy}.

Let $\rho_n$ be the complex spin representation on $\Delta_n$.
The spinor bundle
\[
 \Sigma M=P_{\Spin}(M)\times_{\rho_n}\Delta_n,
 \qquad \rank_{\C}\Sigma M=2^{\lfloor n/2\rfloor},
\]
has the induced Hermitian metric and spin connection $\nabla^\Sigma$.
Clifford multiplication satisfies
\begin{equation}\label{eq:clifford}
 c(X)c(Y)+c(Y)c(X)=-2\ip XY\Id,
 \qquad c(X)^*=-c(X).
\end{equation}
In a local spin frame over an oriented orthonormal frame
$(e_1,\ldots,e_n)$,
\[
 \nabla_X^\Sigma\psi
 =X(\psi)+\frac14\sum_{i,j=1}^n
       \ip{\nabla_Xe_i}{e_j}c(e_i)c(e_j)\psi,
\]
where $X(\psi)$ denotes differentiation of the coefficient functions.
The connection is Hermitian, and its compatibility with Clifford
multiplication takes the form
\[
 \nabla_X^\Sigma(c(Y)\psi)
 =c(\nabla_XY)\psi+c(Y)\nabla_X^\Sigma\psi.
\]
Its curvature is
\begin{equation}\label{eq:spin-curvature}
 R^\Sigma(X,Y)
 =\frac14\sum_{i,j=1}^n
       \ip{R^M(X,Y)e_i}{e_j}c(e_i)c(e_j).
\end{equation}
The intrinsic Dirac operator and connection Laplacian are
\begin{align}
 D_M\psi&=\sum_{i=1}^n c(e_i)\nabla^\Sigma_{e_i}\psi,
       \label{eq:intrinsic-Dirac}\\
 (\nabla^\Sigma)^*\nabla^\Sigma
 &=-\sum_{i=1}^n
   \left(\nabla^\Sigma_{e_i}\nabla^\Sigma_{e_i}
           -\nabla^\Sigma_{\nabla_{e_i}e_i}\right).
       \label{eq:connection-laplacian}
\end{align}
They satisfy the Schr\"odinger--Lichnerowicz formula 
\begin{equation}\label{eq:SL}
 D_M^2=(\nabla^\Sigma)^*\nabla^\Sigma+\frac14\Scal_M.
\end{equation}
Integration over the closed manifold $M$ yields
\begin{equation}\label{eq:integrated-SL}
 \int_M|D_M\psi|^2\vol_M
 =\int_M|\nabla^\Sigma\psi|^2\vol_M
    +\frac14\int_M\Scal_M|\psi|^2\vol_M.
\end{equation}
The self-adjoint operator $D_M$ has discrete real spectrum and smooth
eigenspinors. In general, the square of its first strictly positive
eigenvalue $\lambda_1^+(D_M)$ need not coincide with
$\min\Spec(D_M^2)$, because of a nontrivial kernel or spectral
asymmetry. For a closed spin manifold, the spectrum is symmetric
about zero when $n\not\equiv3\pmod4$, whereas symmetry may fail in
the remaining dimensions \cite[Theorem~1.3.7]{GinouxSpectrum}.
The lower bound proved below applies directly to $D_M^2$. It therefore
controls both signs of the spectrum and excludes the kernel,
independently of spectral symmetry.

In a local orthonormal frame $(e_1,\ldots,e_{m_i})$ of $E_i$, set
\begin{equation}\label{eq:partial-energy}
 |\nabla^\Sigma_{E_i}\psi|^2
   =\sum_{a=1}^{m_i}|\nabla^\Sigma_{e_a}\psi|^2,
 \qquad
 |\nabla^\Sigma\psi|^2
   =\sum_{i=1}^g|\nabla^\Sigma_{E_i}\psi|^2.
\end{equation}
The first expression defines the frame-independent energy density
along $E_i$. The second identity follows from \eqref{eq:splitting}.

To express the relation with ambient spinors, let $\bar c$ denote
spherical Clifford multiplication and use $\bar\nabla$ also for the
spherical spin connection. On the restricted bundle
\begin{equation}\label{restricted-bundle}
 \mathcal E=\Sigma S^{n+1}|_M,
 \qquad N=\bar c(\nu),
 \qquad c_{\mathcal E}(X)=\bar c(X)N,
\end{equation}
the induced hypersurface connection is
\begin{equation}\label{eq:spin-gauss}
 \nabla_X^{\mathcal E}
 =\bar\nabla_X-\frac12\bar c(AX)N.
\end{equation}
Moreover,
\begin{equation}\label{eq:Nparallel}
 N^2=-\Id,\qquad N^*=-N,\qquad \nabla^{\mathcal E}N=0.
\end{equation}
The parallelism assertion follows from $\bar\nabla_XN=-\bar c(AX)$ by \eqref{eq:shape}
and $[-\frac12\bar c(AX)N,N]=\bar c(AX)$.
Define $\mathcal D=\sum_a c_{\mathcal E}(e_a)\nabla^{\mathcal E}_{e_a}$.
Since $N$ is parallel and anticommutes with $c_{\mathcal E}(X)$,
\begin{equation}\label{eq:D-N}
 \mathcal DN=-N\mathcal D.
\end{equation}
Substitution of \eqref{eq:spin-gauss} yields
\begin{equation}\label{eq:extrinsic-D}
 \mathcal D
 =\sum_{a=1}^n\bar c(e_a)N\bar\nabla_{e_a}
  +\frac12\tr A,
\end{equation}
because $\bar c(e_a)N\bar c(e_a)N=-\Id$ in a principal frame.

The standard hypersurface spinor identification provides a
parallel isometric Clifford inclusion
$\iota:\Sigma M\to\mathcal E$ and its orthogonal projection
$\Pi=\iota^*:\mathcal E\to\Sigma M$, satisfying
\begin{equation}\label{eq:intrinsic-projection}
 \Pi\iota=\Id,\qquad \Pi\mathcal D=D_M\Pi.
\end{equation}
Any complementary intrinsic summand carries the opposite Clifford
action. Consequently, $\mathcal D^2$ is unitarily equivalent to a
direct sum of copies of $D_M^2$
\cite[Lemma~2.38 and Theorem~2.39]{BourguignonEtAl}.
The arguments below use these projection and square identities.

\section{Leafwise estimates}\label{sec:leaf-estimates}
Based on the foliations of principal curvature spheres on the isoparametric hypersurface $M^n$ in Section~\ref{sec:conventions}, we derive in this section the leafwise estimates for the spin connection. Let $L$ be a curvature leaf of dimension $m$ with principal curvature
$\kappa$, and set $K=1+\kappa^2=r^{-2}$.
For a local orthonormal frame $(e_1,\ldots,e_m)$ of $TL$, define
the partial Dirac operator on $\Sigma M|_L$ by
\begin{equation}\label{eq:partial-D}
 \DL\psi=\sum_{a=1}^m c(e_a)\nabla^\Sigma_{e_a}\psi.
\end{equation}
Total geodesicity of $L$ implies that the restricted connection satisfies
\[
 \nabla^\Sigma_X(c(Y)\psi)
 =c(\nabla_X^LY)\psi+c(Y)\nabla^\Sigma_X\psi,
 \qquad X,Y\in TL.
\]
Thus $\Sigma M|_L$ is a Hermitian Clifford module with a compatible
connection, and \eqref{eq:partial-D} defines a global operator.
Using metric compatibility and $c(e_a)^*=-c(e_a)$, integration
by parts on the closed leaf yields
\begin{equation}\label{eq:formal-adjoint}
 \int_L\ip{\DL\phi}{\chi}\vol_L
 =\int_L\ip{\phi}{\DL\chi}\vol_L.
\end{equation}

\begin{lemma}\label[lemma]{lem:partial-SL}
On $\Sigma M|_L$,
\begin{equation}\label{eq:leaf-SL}
 \DL^2=(\nabla_L^\Sigma)^*\nabla_L^\Sigma
       +\frac{m(m-1)}4K\Id.
\end{equation}
The adjoint is taken with respect to the Riemannian volume of $L$.
\end{lemma}

\begin{proof}
For $X,Y\in TL$ and $Z\in TM|_L$, substituting $AX=\kappa X$
and $AY=\kappa Y$ into the Gauss equation \eqref{eq:gauss} yields
\begin{equation}\label{eq:R-leaf}
 R^M(X,Y)Z
 =K\bigl(\ip YZ X-\ip XZ Y\bigr).
\end{equation}
In particular, $R^M(X,Y)$ vanishes on the fixed normal space $N_{L/M}$ of $L$ in $M$.

Extend a local orthonormal frame $(e_1,\ldots,e_m)$ of $TL$ to
an orthonormal frame $(e_1,\ldots,e_n)$ of $TM|_L$.
For $a\ne b$, equation \eqref{eq:R-leaf} reads
\[
 R^M(e_a,e_b)e_a=-Ke_b,\qquad
 R^M(e_a,e_b)e_b=Ke_a,
\]
and $R^M(e_a,e_b)e_j=0$ for $j\neq a,b$.
Substitution in \eqref{eq:spin-curvature} yields
\begin{equation}\label{eq:Rspin-leaf}
 \begin{split}
 R^\Sigma(e_a,e_b)
 &=\frac K4\bigl(-c(e_a)c(e_b)+c(e_b)c(e_a)\bigr)\\
 &=-\frac K2c(e_a)c(e_b).
 \end{split}
\end{equation}

For a geodesic tangent frame at a point of $L$, expanding
\eqref{eq:partial-D} in $ \DL^2$ gives
\[
 \DL^2=(\nabla_L^\Sigma)^*\nabla_L^\Sigma
        +\sum_{a<b}c(e_a)c(e_b)R^\Sigma(e_a,e_b).
\]
Since $(c(e_a)c(e_b))^2=-\Id$ when $a\ne b$, the last sum is
\[
 -\frac K2\sum_{a<b}(c(e_a)c(e_b))^2
 =\frac K2\binom m2\Id
 =\frac{m(m-1)}4K\Id.
\]
This proves \eqref{eq:leaf-SL}. The normal curvature term vanishes
by \eqref{eq:R-leaf}.
\end{proof}

\begin{proposition}\label[proposition]{prop:large-leaf}
If $m\geq2$, every smooth section $\psi$ of $\Sigma M|_L$ satisfies
\begin{equation}\label{eq:leaf-bound}
 \int_L|\nabla_L^\Sigma\psi|^2\vol_L
 \geq\frac{mK}{4}\int_L|\psi|^2\vol_L.
\end{equation}
\end{proposition}

\begin{proof}
For each unit vector $e_a$, relation \eqref{eq:clifford} implies
$c(e_a)^*c(e_a)=\Id$. The Cauchy--Schwarz inequality therefore yields
\[
 |\DL\psi|^2
 \leq m\sum_{a=1}^m|c(e_a)\nabla^\Sigma_{e_a}\psi|^2
 =m|\nabla_L^\Sigma\psi|^2.
\]
Formal self-adjointness \eqref{eq:formal-adjoint} and
\cref{lem:partial-SL} imply
\[
 \int_L|\DL\psi|^2\vol_L
 =\int_L|\nabla_L^\Sigma\psi|^2\vol_L
  +\frac{m(m-1)}4K\int_L|\psi|^2\vol_L.
\]
It follows that
\[
 (m-1)\int_L|\nabla_L^\Sigma\psi|^2\vol_L
 \geq\frac{m(m-1)}4K\int_L|\psi|^2\vol_L.
\]
Dividing by $m-1>0$ proves the estimate.
\end{proof}

When $m=1$, the factor $m-1$ vanishes, and the required estimate
depends on the spin holonomy around the curvature circle. A circle
has two spin structures \cite[\S3.3.1]{BourguignonEtAl}. We determine
which one is induced by the embedding.

\begin{lemma}\label[lemma]{lem:circle-holonomy}
Let $L$ be a one-dimensional curvature leaf. Parallel transport for
$\nabla^\Sigma$ once around the embedded circle $L$ acts as $-\Id$
on $\Sigma M|_L$.
\end{lemma}

\begin{proof}
Parametrize $L$ by $x(s)$ with arclength $s\in[0,\ell]$,
$\ell=2\pi r$, and set $T=x'(s)$. Define
\[
 \eta=r(x-\kappa\nu),\qquad
 \zeta=r(\nu+\kappa x).
\]
The vectors $\eta$ and $\zeta$ are orthonormal. As $L$ is a round circle in an affine plane, direct calculations imply
\begin{equation}\label{eq:rotating-frame}
 T'=-\frac{\eta}{r},\qquad
 \eta'=\frac T r,\qquad
 \zeta'=0.
\end{equation}
Here the prime denotes Euclidean differentiation, and $\eta$ is the
outward unit radial vector from the center ($\kappa r\zeta$) of the circle.

Choose the orientation of $L$ and a constant orthonormal basis
$W_1,\ldots,W_{n-1}$ of the fixed normal space $N_{L/M}$, so that the tangent frame
\[
 F(s)=(T(s),W_1,\ldots,W_{n-1})
\]
is positively oriented. Since the leaf is totally geodesic, $\nabla_{T}T=0$, and constancy of $W_{\alpha}$ gives $\nabla_{T}W_{\alpha}=0$. The spin structure \eqref{eq:stable-spin} is obtained by reducing
the Euclidean spin-frame bundle first to the sphere and then to $M$
\cite[\S1]{GinouxProducts}. These reductions are compatible because
each takes the inverse image of an oriented-frame subbundle under
the same double covering.

Consider the stabilized frame
$\widehat F(s)=(T(s),W_1,\ldots,W_{n-1},\nu(s),x(s))$,
with the orientation convention in \eqref{eq:stable-spin}.
The relations
\[
 x=r(\eta+\kappa\zeta),\qquad
 \nu=r(-\kappa\eta+\zeta)
\]
relate $\widehat F$ by a constant change of basis and a fixed
reordering to
\[
 (T,\eta,\zeta,W_1,\ldots,W_{n-1}).
\]
Equation \eqref{eq:rotating-frame} shows that this frame makes one
full $2\pi$ rotation in the fixed plane spanned by $T$ and $\eta$,
while all other vectors remain constant. Embeddedness ensures that
the parametrization of length $2\pi r$ traverses the circle exactly once.

A plane rotation through angle $\theta$ lifts to the spin group as
\[
 \cos(\theta/2)+\sin(\theta/2)\,u v,
\]
where $u,v$ form an orthonormal basis of the plane and the product is
taken in the real Clifford algebra. At $\theta=2\pi$, this lift equals
$-1$. The lift of the loop $\widehat F$ therefore ends at minus its
initial point.

Under \eqref{eq:stable-spin}, the lift of $F$ to the spin-frame
bundle of $M$ has the same endpoint sign. The stabilized frame
$\widehat F$ rotates in Euclidean space, whereas $F$ is parallel under the intrinsic
connection as shown above. Its spin lift is therefore horizontal, and the endpoint
of this lift determines the intrinsic holonomy. The central element
$-1$ acts by $-\Id$ on complex spinors, so the holonomy is $-\Id$.
\end{proof}

\begin{proposition}\label[proposition]{prop:circle-bound}
If $m=1$, every smooth section $\psi$ of $\Sigma M|_L$ satisfies
\[
 \int_L|\nabla_T^\Sigma\psi|^2\vol_L
 \geq\frac1{4r^2}\int_L|\psi|^2\vol_L.
\]
Thus \eqref{eq:leaf-bound} also holds for one-dimensional leaves.
\end{proposition}

\begin{proof}
Choose a parallel orthonormal spinor frame over $[0,\ell]$ (as above, $L$ is parameterized by $x(s)$, $s\in [0,\ell]$, $\ell=2\pi r$).
By \cref{lem:circle-holonomy}, the coefficient vector
$f:[0,\ell]\to\C^{\rank\Sigma M}$ of $\psi$ satisfies
$f(\ell)=-f(0)$. Its derivatives satisfy the same endpoint
condition. Thus $f$ has a smooth antiperiodic extension to $\R$,
with Fourier series
\[
 f(s)=\sum_{k\in\Z}a_k
      \exp\!\left(\frac{2\pi\mathrm i(k+\tfrac12)s}{\ell}\right).
\]
Parseval's identities imply
\begin{align*}
 \int_0^\ell |f(s)|^2\dd s
   &=\ell\sum_{k\in\Z}|a_k|^2,\\
 \int_0^\ell |f'(s)|^2\dd s
   &=\ell\sum_{k\in\Z}
        \left(\frac{2\pi(k+\tfrac12)}{\ell}\right)^2|a_k|^2\\
   &\geq\left(\frac{\pi}{\ell}\right)^2
        \int_0^\ell|f(s)|^2\dd s.
\end{align*}
The frame is parallel and orthonormal, so $|f|=|\psi|$ and
$|f'|=|\nabla_T^\Sigma\psi|$. Substituting $\ell=2\pi r$ shows
that the coefficient is $1/(4r^2)$, which proves the assertion.
\end{proof}

Combining \cref{prop:large-leaf,prop:circle-bound}, we obtain,
for every principal curvature $\kappa_i$ with multiplicity $m_i$,
\begin{equation}\label{eq:leaf-all}
 \int_L|\nabla_L^\Sigma\psi|^2\vol_L
 \geq\frac{m_i(1+\kappa_i^2)}4
       \int_L|\psi|^2\vol_L
 \quad\text{for every curvature leaf }L\text{ of }E_i.
\end{equation}
This estimate applies to every smooth spinor restricted to a leaf,
without any orthogonality condition. It can therefore be integrated
over $M$ and summed over the principal distributions.

\section{The connection estimate}\label{sec:global}

We integrate the leafwise estimates by applying the coarea formula
to the focal maps. Their Jacobians are constant, so the resulting
inequalities retain the leafwise constants. Throughout this section,
Hausdorff measure is normalized to agree with the induced Riemannian
volume on smooth Euclidean submanifolds.

\begin{lemma}\label[lemma]{lem:coarea}
Let $F:P^n\to\R^N$ be a smooth map of constant rank $q$, where
$P$ is a closed Riemannian manifold. For $q>0$, let $J_qF$ denote
the product of the nonzero singular values of $dF$ under the metrics of $P^n$ and $\mathbb{R}^{N}$. Then every nonnegative Borel function $h$ satisfies
\begin{equation}\label{eq:coarea}
 \int_P hJ_qF\vol_P
 =\int_{F(P)}
     \left(\int_{F^{-1}(y)}h\vol_{F^{-1}(y)}\right)
     d\HS^q(y).
\end{equation}
Each fiber is an embedded submanifold of dimension $n-q$.
For $q=0$, the same formula holds with $J_0F=1$ and counting
measure $\HS^0$ on the finite image.
\end{lemma}

\begin{proof}
Suppose $q>0$. The constant-rank theorem and compactness provide
a finite cover of $P$ by coordinate neighborhoods $U_\alpha$ such
that $F|_{U_\alpha}:U_\alpha\to N_\alpha$ is a submersion onto an
embedded $q$-dimensional submanifold $N_\alpha=F(U_\alpha)\subset\R^N$.
Choose a subordinate nonnegative smooth partition of unity
$(\rho_\alpha)$.

In constant-rank coordinates $(u,v)$, write $F(u,v)=z(u)$ and
represent the domain metric by a block matrix $G=\left(\begin{array}{cc}
                                                   G_{uu} & G_{uv} \\
                                                   G_{vu} & G_{vv} 
                                                 \end{array}\right)$. 
Let $B$ be the matrix of the induced Euclidean metric on $N_\alpha$ in the
$u$ coordinates. The fiber metric is $G_{vv}$.
Set $Q=G_{uu}-G_{uv}G_{vv}^{-1}G_{vu}$. For $q=n$, use the
conventions $\det G_{vv}=1$ and $Q=G_{uu}=G$. The block determinant
identities yield
\[
 \det G=\det G_{vv}\det Q,\qquad
 (J_qF)^2=\det B\,\det((G^{-1})_{uu})
          =\frac{\det B}{\det Q}.
\]
It follows that
$J_qF\sqrt{\det G}=\sqrt{\det B}\sqrt{\det G_{vv}}$.
Tonelli's theorem applied to $\rho_\alpha h$ then gives
\[
 \int_{U_\alpha}\rho_\alpha hJ_qF\vol_P
 =\int_{N_\alpha}
   \left(\int_{F^{-1}(y)\cap U_\alpha}
        \rho_\alpha h\vol_{F^{-1}(y)}\right)d\HS^q(y).
\]
The Jacobian in this identity is $J_qF$, because the target tangent
space is the image of $dF$ equipped with its induced metric. Extend each
fiber integral by zero outside $N_\alpha$ and sum over $\alpha$.
On each full fiber, the integrand becomes
$h\sum_\alpha\rho_\alpha=h$. This proves \eqref{eq:coarea},
including points where image patches overlap. The argument therefore
also applies when the full image $F(P)$ is not embedded.

If $q=0$, then $F$ is locally constant and has finite image by
compactness. In this case, the formula reduces to a sum over the fibers.
\end{proof}

\begin{lemma}\label[lemma]{lem:focal}
For each principal curvature $\kappa_i$, define the focal map
\begin{equation}\label{eq:focal}
 F_i=\frac{\nu+\kappa_i x}{\sqrt{1+\kappa_i^2}}
       :M\longrightarrow S^{n+1}(1)\subset\R^{n+2}.
\end{equation}
Then $\ker dF_i=E_i$, $\rank dF_i=n-m_i$, and
\begin{equation}\label{eq:focal-J}
 J_i:=J_{n-m_i}F_i
 =\prod_{j\ne i}
   \left(\frac{|\kappa_i-\kappa_j|}
               {\sqrt{1+\kappa_i^2}}\right)^{m_j}.
\end{equation}
With the empty product interpreted as $1$, the Jacobian $J_i$ is
positive and constant. Each connected component of a fiber is a
curvature leaf of $E_i$.
\end{lemma}

\begin{proof}
The vector $\nu+\kappa_i x$ has squared norm $1+\kappa_i^2$.
As $\kappa_i$ is constant, differentiation yields
\[
 dF_i(X)=\frac{(\kappa_i\Id-A)X}{\sqrt{1+\kappa_i^2}}.
\]
On $E_j$, this is multiplication by
$(\kappa_i-\kappa_j)/\sqrt{1+\kappa_i^2}$.
The distributions $E_j$ are mutually orthogonal. Hence the nonzero
absolute values of these factors are the singular values, with
multiplicities $m_j$. The formulas for the kernel, rank, and Jacobian
follow.

By the constant-rank theorem, each fiber is a smooth embedded
$m_i$-dimensional submanifold tangent to $E_i$. It is compact,
being closed in $M$, and each of its connected components is an
integral submanifold of $E_i$. Conversely, $dF_i$ vanishes on $E_i$,
so every connected curvature leaf lies in a single fiber. Thus the
maximal leaves are exactly the connected components of the fibers.
These components are open in the fiber. A finite cover by connected
coordinate neighborhoods exists by compactness, so each fiber has
only finitely many components.
\end{proof}

\begin{proposition}\label[proposition]{prop:partial-global}
For every $i$ and every $\psi\in\Gamma(\Sigma M)$,
\begin{equation}\label{eq:partial-global}
 \int_M|\nabla^\Sigma_{E_i}\psi|^2\vol_M
 \geq\frac{m_i(1+\kappa_i^2)}4\int_M|\psi|^2\vol_M.
\end{equation}
\end{proposition}

\begin{proof}
By \cref{lem:focal}, each fiber of $F_i$ is a finite disjoint union
of curvature leaves. Summing \eqref{eq:leaf-all} over these components,
we obtain, for $y\in F_i(M)$,
\[
 \int_{F_i^{-1}(y)}|\nabla^\Sigma_{E_i}\psi|^2\vol
 \geq \frac{m_i(1+\kappa_i^2)}4\int_{F_i^{-1}(y)}|\psi|^2\vol.
\]
Applying \cref{lem:coarea} to the two nonnegative integrands yields
\[
 J_i\int_M|\nabla^\Sigma_{E_i}\psi|^2\vol_M
 \geq \frac{m_i(1+\kappa_i^2)}4J_i\int_M|\psi|^2\vol_M.
\]
Division by the positive constant $J_i$ proves the estimate.
When $n-m_i=0$, we have $E_i=TM$, and connectedness forces $F_i$
to be constant. The assertion then reduces to \eqref{eq:leaf-all}
on the leaf $M$.
\end{proof}

\begin{proof}[Proof of the inequality in \cref{thm:connection}]
Sum \eqref{eq:partial-global} over the principal distributions.
Using \eqref{eq:partial-energy} and the identity
$\sum_i m_i(1+\kappa_i^2)=n+|A|^2$, we obtain
\[
 \int_M|\nabla^\Sigma\psi|^2\vol_M
 \geq \frac14\sum_i m_i(1+\kappa_i^2)\int_M|\psi|^2\vol_M
 =\frac{n+|A|^2}{4}\int_M|\psi|^2\vol_M.
\]
The endpoint eigenspinors constructed in \cref{sec:killing} will
attain this bound and establish its optimality.
\end{proof}

\begin{proposition}\label[proposition]{prop:spectral-gap}
For every closed isoparametric hypersurface
with the induced spin structure,
\begin{equation}\label{eq:global-gap}
 \int_M|D_M\psi|^2\vol_M
 \geq \frac{n^2}{4}(1+H^2)\int_M|\psi|^2\vol_M
 \qquad(\psi\in\Gamma(\Sigma M)).
\end{equation}
\end{proposition}

\begin{proof}
Substituting \eqref{eq:scalar} into \eqref{eq:integrated-SL} gives
\[
 \int_M|D_M\psi|^2\vol_M
 =\int_M|\nabla^\Sigma\psi|^2\vol_M
  +\frac{n(n-1)+n^2H^2-|A|^2}{4}\int_M|\psi|^2\vol_M.
\]
After applying \eqref{eq:connection-main}, the coefficient of
$\int_M|\psi|^2\vol_M$ in the lower bound becomes
\[
 \frac{n+|A|^2+n(n-1)+n^2H^2-|A|^2}{4}
 =\frac{n^2}{4}(1+H^2).
\]
This proves \eqref{eq:global-gap}.
\end{proof}
 For a nonzero eigenspinor satisfying
$D_M\psi=\lambda\psi$, division by $\int_M|\psi|^2$ in \eqref{eq:global-gap} yields
$\lambda^2\geq n^2(1+H^2)/4$. Consequently, the kernel is trivial
and the spectrum is disjoint from
$(-\frac n2\sqrt{1+H^2},\frac n2\sqrt{1+H^2})$. This, together with the endpoint eigenspinors constructed in \cref{sec:killing}, proves \cref{cor:cmc} and \cref{thm:minimal}.

\section{Endpoint eigenvalues and the Dirac spectrum}\label{sec:killing}

We now construct eigenspinors at both endpoints of the spectral gap.
The construction applies to closed connected embedded hypersurfaces
of constant mean curvature in the unit sphere, oriented by a chosen
unit normal. It uses the evaluation isomorphism for ambient Killing
spinors together with the projection \eqref{eq:intrinsic-projection}.
We first recall this evaluation property
\cite[\S2, Lemma 1]{BaerSphere}.

\begin{lemma}\label[lemma]{lem:killing-evaluation}
On $S^{n+1}(1)$, $n\geq1$, evaluation at any point is an isomorphism
from the space of Killing spinors with Killing constant $1/2$ onto
the spinor fiber at that point. In particular, for every $p$ and
$v\in\Sigma_pS^{n+1}$, there is a unique smooth spinor $\Psi$
satisfying
\begin{equation}\label{eq:killing}
 \bar\nabla_X\Psi=\frac12\bar c(X)\Psi
 \quad\text{for all }X,\qquad \Psi(p)=v.
\end{equation}
\end{lemma}

\begin{proof}
Consider the connection
$\widehat\nabla_X=\bar\nabla_X-\frac12\bar c(X)$ on the ambient spinor
bundle. Torsion-freeness and Clifford compatibility imply that its
curvature is
\[
 \widehat R(X,Y)
 =\bar R^\Sigma(X,Y)
  +\frac14[\bar c(X),\bar c(Y)].
\]
On the unit sphere, formula \eqref{eq:spin-curvature} reads
$\bar R^\Sigma(X,Y)=-\frac14[\bar c(X),\bar c(Y)]$, and hence
$\widehat R=0$. The sphere $S^{n+1}$ is simply connected for $n\geq1$,
so parallel transport is path-independent. It therefore extends $v$
uniquely to a global $\widehat\nabla$-parallel section, which proves
\eqref{eq:killing}. Clifford multiplication is skew-Hermitian, so
$\widehat\nabla$ is Hermitian. A nonzero initial value thus determines
a nowhere-vanishing spinor.
\end{proof}

We use the restricted bundle $\mathcal E$ and the operators $N$
and $\mathcal D$ defined in \cref{sec:conventions}.

\begin{proposition}\label[proposition]{prop:endpoints}
Let $M^n\subset S^{n+1}(1)$ be a closed connected embedded
hypersurface. Choose a unit normal, orient $M$ accordingly, and
equip it with the induced spin structure. If the normalized mean
curvature $H$ is constant, then
\[
 \left\{-\frac n2\sqrt{1+H^2},
          \frac n2\sqrt{1+H^2}\right\}
 \subset\Spec(D_M).
\]
\end{proposition}

\begin{proof}
Set $a=nH/2$, $b=n/2$, and $\Lambda=\sqrt{a^2+b^2}=\frac n2\sqrt{1+H^2}$.
Let $\Psi$ be an ambient Killing spinor as in
\cref{lem:killing-evaluation}, and write $\phi=\Psi|_M$.
Substituting $\bar c(e_a)N\bar c(e_a)=N$ into
\eqref{eq:extrinsic-D}, we obtain
\begin{equation}\label{eq:two-by-two}
 \mathcal D\phi=a\phi+bN\phi,\qquad
 \mathcal D(N\phi)=b\phi-aN\phi,
\end{equation}
where the second identity follows from \eqref{eq:D-N}. Both coefficients
$a$ and $b$ are constant.

For $\lambda=\Lambda$ or $\lambda=-\Lambda$, set
\[
 \phi_\lambda=(\lambda+a)\phi+bN\phi.
\]
Using \eqref{eq:two-by-two} and $\lambda^2=a^2+b^2$,
\begin{align*}
 \mathcal D\phi_\lambda
 &=\bigl(a(\lambda+a)+b^2\bigr)\phi
     +\bigl(b(\lambda+a)-ab\bigr)N\phi\\
 &=\lambda(\lambda+a)\phi+\lambda bN\phi
 =\lambda\phi_\lambda.
\end{align*}
Moreover, $b>0$ implies $\Lambda>|a|$, so $\lambda+a\ne0$.

The endomorphism $T_\lambda=(\lambda+a)\Id+bN$ is invertible,
with inverse
\[
 T_\lambda^{-1}
 =\frac{(\lambda+a)\Id-bN}{(\lambda+a)^2+b^2}.
\]
Fix $p\in M$ and choose $w\in\mathcal E_p$ with $\Pi w\ne0$.
The evaluation isomorphism in \cref{lem:killing-evaluation} allows
us to choose $\Psi$ with $\Psi(p)=T_\lambda^{-1}w$.
Equation \eqref{eq:intrinsic-projection} then yields
\[
 D_M(\Pi\phi_\lambda)=\lambda\Pi\phi_\lambda,
 \qquad (\Pi\phi_\lambda)(p)=\Pi w\ne0.
\]
Thus each endpoint is an eigenvalue of the intrinsic Dirac operator.
The argument applies in every dimension.
\end{proof}

\begin{proof}[Completion of the proofs of
\cref{thm:minimal,thm:connection,cor:cmc}]
For an isoparametric hypersurface, \eqref{eq:splitting} shows that
$H$ is constant. Proposition~\ref{prop:spectral-gap} excludes all
eigenvalues in $(-\Lambda,\Lambda)$, while \cref{prop:endpoints}
realizes both endpoints. This proves \cref{cor:cmc}. Setting $H=0$
proves \cref{thm:minimal}.

To prove optimality in \cref{thm:connection}, choose a nonzero
endpoint eigenspinor $\psi$ from \cref{prop:endpoints}. The scalar
curvature $\Scal_M$ is constant, so \eqref{eq:SL} implies
\[
 (\nabla^\Sigma)^*\nabla^\Sigma\psi
 =\left(\Lambda^2-\frac{\Scal_M}{4}\right)\psi
 =\frac{n+|A|^2}{4}\psi.
\]
The spinor $\psi$ attains equality in \eqref{eq:connection-main}.
Thus the constant is optimal, and
\[
 \min\Spec\bigl((\nabla^\Sigma)^*\nabla^\Sigma\bigr)
 =\frac{n+|A|^2}{4}.
\]
\end{proof}

\section{From the Dirac spectral gap to the first Laplace eigenvalue}
\label{sec:laplace}

Let $x:M^n\hookrightarrow S^{n+1}(1)$ be the inclusion of a closed
connected embedded minimal hypersurface, with $n\ge1$. Throughout
this section, we use the normal and spinor conventions of
\cref{sec:conventions}, but impose no condition on the principal
curvatures beyond minimality. Write $\Delta=\Delta_M$.
All differential forms in this section are real. The notation
$\norm{\cdot}$ and $\Ltwo{\cdot}{\cdot}$ refers to the $L^2$ norm
and real $L^2$ inner product, respectively. Minimality and
Takahashi's formula yield
\begin{equation}\label{lap:coordinate}
 \tr A=0,\qquad \Delta x=nx.
\end{equation}

To relate the Dirac and Laplace operators, we need a faithful representation of the complex
Clifford algebra and consider an isometric embedding of the exterior algebra of differential forms into this representation (\cref{lap:lem:trace}).
For odd $n$, the restricted bundle $\mathcal E\simeq\Sigma M\oplus\Sigma M$
already carries the faithful representation, so one copy suffices for the isometry. For even $n$, it is no longer faithful.
Hence we use the doubled bundle below to give a unified proof of
\cref{lap:lem:trace} for both odd and even dimensions.

To distinguish the two form components in the Clifford trace, we use
two copies of the restricted ambient spinor bundle $\mathcal E$. We
keep the same unit normal $\nu$ on both copies and take the ambient
Clifford action to be $\widehat{\bar c}=\bar c\oplus(-\bar c)$.
Since $N=\bar c(\nu)$ and $c_{\mathcal E}(X)=\bar c(X)N$ as in \eqref{restricted-bundle}, this gives the doubled bundle
\[
 \E=\mathcal E\oplus\mathcal E,
 \quad c(X)=c_{\mathcal E}(X)\oplus c_{\mathcal E}(X),
 \quad \hat{N}=N\oplus(-N), \quad c(X)=\widehat{\bar c}(X)\hat{N}.
\]

Write $\bar\nabla^{\oplus}=\bar\nabla\oplus\bar\nabla$ for the
ambient spin connection on the doubled bundle. Applying the spinorial Gauss
formula \eqref{eq:spin-gauss} on each summand gives the induced
connection $\nabla=\nabla^{\mathcal E}\oplus\nabla^{\mathcal E}$ on $\E$ (doubled Gauss formula)
\begin{equation}\label{eq-doubleGauss}
 \nabla_X
 =\bar\nabla_X^{\oplus}-\frac12\widehat{\bar c}(AX)\hat{N}
 =\bar\nabla_X^{\oplus}-\frac12c(AX).
\end{equation}
In particular,
\[
 \nabla \hat{N}=0,\qquad \hat{N}^2=-\Id,\qquad \hat{N}c(X)=-c(X)\hat{N},
 \qquad \De=\mathcal D\oplus\mathcal D.
\]
The square identity following \eqref{eq:intrinsic-projection}
transfers the assumed gap \eqref{eq:transfer-gap} to
\begin{equation}\label{lap:quadratic-gap}
 \norm{\De\psi}^2\ge \frac{n^2}{4}\norm{\psi}^2
 \qquad\text{for every smooth section }\psi\text{ of }\E.
\end{equation}
Put $r=\rank_{\C}\E$. For a differential form $\alpha$, let
$c(\alpha)$ denote its Clifford action: an increasing wedge
monomial is sent to the corresponding Clifford product.

\begin{lemma}\label[lemma]{lap:lem:trace}
For each $q\in M$, the real-linear map
\[
 \begin{aligned}
 \Lambda^*T_q^*M\oplus\Lambda^*T_q^*M
 &\longrightarrow \operatorname{End}_{\C}(\E_q),\\
 (\alpha,\beta)&\longmapsto c(\alpha)+c(\beta)\hat{N}
 \end{aligned}
\]
is an isometric embedding for the normalized real Hilbert--Schmidt
inner product on the target. More precisely,
\begin{equation}\label{lap:eq:trace-isometry}
 \frac1r\Rea\tr\bigl((c(\alpha)+c(\beta)\hat{N})^*
                         (c(\alpha)+c(\beta)\hat{N})\bigr)
 =|\alpha|^2+|\beta|^2.
\end{equation}
The corresponding inner product identity follows by polarization.
\end{lemma}

\begin{proof}
The generators $c(e_1),\ldots,c(e_n),\hat{N}$ satisfy the Clifford
relations in dimension $n+1$. Let $m_I$ denote an ordered product
of distinct generators, including the empty product which is the identity. Each $m_I$
is unitary, and for $I\ne J$ the product $m_I^*m_J$ is, up to
sign, a nonempty Clifford monomial.

A nonempty tangential monomial of odd degree anticommutes with
the normal multiplication on each copy of $\mathcal E$. One of
even degree anticommutes with a tangent generator occurring in
it. In either case its trace vanishes. A monomial containing $\hat{N}$
has opposite traces on the two copies of $\mathcal E$, because
the tangent actions agree and the normal actions have opposite
signs. Its total trace therefore vanishes as well. The monomials
$m_I$ are consequently orthonormal for the normalized real trace
pairing, which proves \eqref{lap:eq:trace-isometry} and injectivity.
\end{proof}

\begin{lemma}\label[lemma]{lap:lem:killing-frame}
There is a global unitary frame $s_1,\ldots,s_r$ of
$\E$ satisfying
\begin{equation}\label{lap:eq:frame-equation}
 \nabla_Xs_j=-\frac12c(X)s_j+\frac12c(AX)\hat{N}s_j,
 \qquad \De s_j=\frac{n}{2}s_j.
\end{equation}
\end{lemma}

\begin{proof}
On the sphere, the connection
$\bar\nabla^{\oplus}-\widehat{\bar c}/2$ is metric and flat: the proof
of \cref{lem:killing-evaluation} applies to either sign of the
ambient Clifford action. It therefore admits a global unitary
Killing spinor frame $\rho_1,\ldots,\rho_r$
satisfying
\[
 \bar\nabla_X^{\oplus}\rho_j=\frac12\widehat{\bar c}(X)\rho_j.
\]
This is the standard construction of spherical Killing spinors
\cite{BaerSphere,LawsonMichelsohn}. Restrict the frame to $M$. Since
$\widehat{\bar c}(X)=-c(X)\hat{N}$, the doubled Gauss formula \eqref{eq-doubleGauss} gives
\[
 \nabla_X\rho_j
 =\frac12\widehat{\bar c}(X)\rho_j-\frac12c(AX)\rho_j
 =-\frac12\bigl(c(X)\hat{N}+c(AX)\bigr)\rho_j.
\]
Now set $U=(\Id+\hat{N})/\sqrt2$ and $s_j=U\rho_j$. The identities
$\nabla\hat{N}=0$, $\hat{N}^*=-\hat{N}$, and $\hat{N}^2=-\Id$ show
that $U$ is parallel and unitary, with $U^{-1}=(\Id-\hat{N})/\sqrt2$.
For $Y\in TM$, anticommutation with $\hat{N}$ gives
\[
 Uc(Y)U^{-1}
 =\tfrac12c(Y)(\Id-\hat{N})^2=-c(Y)\hat{N},
 \qquad U\bigl(c(Y)\hat{N}\bigr)U^{-1}=c(Y).
\]
Conjugating the two terms in $\nabla_X\rho_j$ now yields
\begin{align*}
 \nabla_Xs_j
 &=-\tfrac12 U\bigl(c(X)\hat{N}+c(AX)\bigr)U^{-1}s_j\\
 &=-\tfrac12c(X)s_j+\tfrac12c(AX)\hat{N}s_j.
\end{align*}
It follows that
\[
 \De s_j=\frac n2s_j
 +\frac12\sum_i c(e_i)c(Ae_i)\hat{N}s_j
 =\frac n2s_j,
\]
because $\sum_i c(e_i)c(Ae_i)=-\tr A=0$.
\end{proof}

We henceforth write Clifford multiplication by juxtaposition and
identify tangent vectors with their metric dual one-forms.
Applying \cref{lap:lem:trace} in the unitary frame of
\cref{lap:lem:killing-frame}, we obtain
\begin{equation}\label{lap:eq:frame-trace}
 \frac1r\sum_{j=1}^r\norm{(\alpha+\beta \hat{N})s_j}^2
 =\norm{\alpha}^2+\norm{\beta}^2.
\end{equation}
The doubled bundle thus provides a faithful complex Clifford
representation and the same trace identity in every dimension.

The shape operator acts on $p$-forms by the symmetric endomorphism
$A^{[p]}$, defined for $0\le p\le n$ by
\begin{equation}\label{lap:eq:A-extension}
 (A^{[p]}\eta)(X_1,\ldots,X_p)
 =\sum_{j=1}^p\eta(X_1,\ldots,AX_j,\ldots,X_p).
\end{equation}
If $Ae_i=\mu_i e_i$, for $i=1,\ldots,n$, then in the dual coframe $\{e^i\}$, the action on the basis
$\bigl\{e^I = e^{i_1}\wedge\cdots\wedge e^{i_p}
\bigm| I=(i_1,\ldots,i_p),\ 1\leq i_1<\cdots<i_p\leq n\bigr\}$
of $p$-forms is
\[
 A^{[p]}e^I=\Big(\sum_{i\in I}\mu_i\Big)e^I,
 \qquad |I|=p.
\]
The extension in degree zero vanishes by definition, while its
action on top-degree forms is multiplication by $\tr A=0$. Thus
\begin{equation}\label{lap:eq:endpoint-A}
 A^{[0]}=0,\qquad A^{[n]}=0.
\end{equation}
Let $\delta$ be the formal adjoint of $d$. On
$\Omega^*(M)=\bigoplus_{p=0}^n\Omega^p(M)$, define
\begin{equation}\label{lap:eq:operators}
 \KK|_{\Omega^p}=(-1)^pA^{[p]},\qquad
 \NN|_{\Omega^p}=(-1)^p\bigl(\tfrac{n}{2}-p\bigr).
\end{equation}

\begin{proposition}\label[proposition]{lap:prop:block}
For smooth differential forms $\alpha$ and $\beta$,
\begin{equation}\label{lap:eq:block-action}
 \begin{split}
 \De\bigl((\alpha+\beta \hat{N})s_j\bigr)
 =\bigl[&(d+\delta)\alpha+\NN\alpha+\KK\beta\\
        &+((d+\delta)\beta-\NN\beta+\KK\alpha)\hat{N}\bigr]s_j.
 \end{split}
\end{equation}
The corresponding operator on the form coefficients is
\begin{equation}\label{lap:eq:block-matrix}
(\alpha,\beta)\mapsto \begin{pmatrix}
 d+\delta+\NN&\KK\\
 \KK&d+\delta-\NN
 \end{pmatrix}\begin{pmatrix}\alpha\\\beta\end{pmatrix}.
\end{equation}
\end{proposition}

\begin{proof}
Fix a point and choose a local orthonormal frame whose covariant
derivative vanishes there. Recall \cite{LawsonMichelsohn} that $d=\sum_i e^i\wedge\nabla_{e_i}$ and $\delta=-\sum_i\iota_{e_i}\nabla_{e_i}$. With respect to the canonical isomorphism $Cl_n\cong \Lambda^*(\mathbb{R}^n)$ as vector spaces, Clifford multiplication between $v\in \mathbb{R}^n$ and $\varphi\in Cl_n$ can be written as $$v\varphi=v\wedge \varphi-\iota_{v}\varphi.$$ Therefore, we have
\[
 \sum_{i=1}^n e_i\,\nabla_{e_i}\alpha
 =\sum_{i=1}^n\bigl(e^i\wedge\nabla_{e_i}\alpha
                  -\iota_{e_i}\nabla_{e_i}\alpha\bigr)
 =d\alpha+\delta\alpha.
\]
Applying
\eqref{lap:eq:frame-equation}, we find
\[
 \De(\alpha s_j)
 =\left((d+\delta)\alpha-\frac12\sum_i e_i\alpha e_i
        +\frac12\sum_i e_i\alpha Ae_i \hat{N}\right)s_j.
\]
For a $p$-form $\eta$, the Clifford relations imply
\begin{equation}\label{lap:eq:sandwich}
 \sum_i e_i\eta e_i=(-1)^p(2p-n)\eta=-2\NN\eta,
 \quad
 \sum_i e_i\eta Ae_i=2(-1)^pA^{[p]}\eta=2\KK\eta.
\end{equation}
To verify the second identity, diagonalize $A$. On a monomial
$e^I$, the terms with $i\in I$ have sign $(-1)^p$, whereas those
with $i\notin I$ have sign $(-1)^{p+1}$. The identity follows
from $\sum_i\mu_i=0$. Substituting both identities proves the
formula for $\De(\alpha s_j)$.

For $\beta \hat{N}$, parallelism of $\hat{N}$ leaves the derivative term
$((d+\delta)\beta)\hat{N}$. Anticommuting $\hat{N}$ with $e_i$ reverses the
sign of the $\NN$ term, whereas $\hat{N}(Ae_i)\hat{N}=Ae_i$ preserves the
sign of the $\KK$ term. These observations complete
\eqref{lap:eq:block-action}.
\end{proof}

Set $\LA=d+\delta+\KK$. On the closed manifold $M$, the operator
$d+\delta$ is self-adjoint with domain $H^1$ and has compact
resolvent. The smooth symmetric zeroth-order endomorphism $\KK$
is bounded on $L^2$. Hence $\LA$ is self-adjoint on the same
domain and has compact resolvent. Its principal symbol is that of
$d+\delta$, so it is elliptic.

\begin{proposition}\label[proposition]{lap:prop:quadratic}
Let $\chi=\sum_{p=0}^n\chi_p \in \Omega^{*}(M)=\bigoplus_{p=0}^{n}\Omega^{p}(M)$, and denote its odd and even parts
by $\alpha$ and $\beta$, respectively. For a smooth real function
$u$, set
\[
 \psi_j=(u+\alpha+\beta \hat{N})s_j.
\]
Then
\begin{equation}\label{lap:eq:quadratic-identity}
 \begin{split}
 \frac1r\sum_{j=1}^r
 \bigl(\norm{\De\psi_j}^2-\frac{n^2}{4}\norm{\psi_j}^2\bigr)
 ={}&\norm{\LA\chi}^2
      -\sum_{p=0}^n p(n-p)\norm{\chi_p}^2\\
    &+\norm{du}^2+2\Ltwo{du}{\chi_1}.
 \end{split}
\end{equation}
\end{proposition}

\begin{proof}
Take $u=0$ first. In \eqref{lap:eq:block-action},
$(d+\delta)\alpha+\KK\beta$ and $(d+\delta)\beta+\KK\alpha$
are the even and odd parts of $\LA\chi$. In their respective
components, they are orthogonal to $\NN\alpha$ and $\NN\beta$
by parity and \cref{lap:lem:trace}. The identity
\eqref{lap:eq:frame-trace} therefore yields
\[
 \frac1r\sum_j\norm{\De\psi_j}^2
 =\norm{\LA\chi}^2+
   \sum_{p=0}^n\bigl(\tfrac{n}{2}-p\bigr)^2\norm{\chi_p}^2.
\]
Subtracting $\tfrac{n^2}{4}\norm{\chi}^2=\frac1r\sum_j\tfrac{n^2}{4}\norm{\psi_j}^2$ gives the first line of
\eqref{lap:eq:quadratic-identity}.

For general $u$, put $\phi_j=(\alpha+\beta\hat{N})s_j$. By
\eqref{lap:eq:frame-equation},
\[
 \psi_j=\phi_j+us_j,\qquad
 \De(us_j)=\bigl(du+\tfrac n2u\bigr)s_j.
\]
In the first component of $\De\phi_j$ in \eqref{lap:eq:block-action},
the terms of degrees zero and one are $\delta\chi_1$ and
$(1-\tfrac n2)\chi_1$, respectively. Expanding the squares by
\eqref{lap:eq:frame-trace} gives
\begin{align*}
 \frac1r\sum_j\bigl(\norm{\psi_j}^2-\norm{\phi_j}^2\bigr)
 &=\norm{u}^2,\\
 \frac1r\sum_j\bigl(\norm{\De\psi_j}^2-\norm{\De\phi_j}^2\bigr)
 &=\tfrac{n^2}{4}\norm{u}^2+\norm{du}^2\\
 &\quad+n\Ltwo{u}{\delta\chi_1}+(2-n)\Ltwo{du}{\chi_1}\\
 &=\tfrac{n^2}{4}\norm{u}^2+\norm{du}^2+2\Ltwo{du}{\chi_1},
\end{align*}
where $\Ltwo{u}{\delta\chi_1}=\Ltwo{du}{\chi_1}$. Subtracting
$n^2/4$ times the first identity from the second and using the case
$u=0$ proves \eqref{lap:eq:quadratic-identity}.
\end{proof}

On the closed manifold $M$, the Hodge decomposition of smooth real
$p$-forms is orthogonal for the $L^2$ inner product, giving
\[
 \Omega^p(M)=d\Omega^{p-1}(M)\oplus
             \delta\Omega^{p+1}(M)\oplus\mathcal H^p(M),
 \qquad \mathcal H^p(M)=\ker d\cap\ker\delta,
\]
where terms of degree below zero or above $n$ are omitted. 

In degree zero, connectedness gives $\mathcal H^0(M)=\R\cdot1$. Thus
\[
 C^\infty(M)=\R\cdot1\oplus C^\infty_\perp(M),\qquad
 C^\infty_\perp(M):=\left\{f\in C^\infty(M):\int_M f\,\volM=0\right\}.
\]
On $C^\infty(M)$ the Laplacian $\Delta=\delta d+d\delta$ has kernel $\R\cdot1$ and
image $\delta\Omega^1(M)=C^\infty_\perp(M)$, the space of mean-zero functions.
It restricts to a bijection of $C^\infty_\perp(M)$ onto itself,
whose inverse we denote by $\Delta^{-1}$.

In $L^2(T^*M)$, the closure of the exact one-forms is $dH^1(M)$ by
the Poincar\'e inequality on mean-zero functions. We denote the
orthogonal projection onto this subspace by $\Pe$.
For a smooth one-form $\eta$, we have
$\Ltwo{\delta\eta}{1}=\Ltwo{\eta}{d1}=0$. Write its Hodge
decomposition as $\eta=du+\delta\gamma+h$, with $h\in\mathcal H^1(M)$,
and choose $u$ to have mean zero. Since $\delta^2=0$ and $\delta h=0$,
applying $\delta$ gives $\Delta u=\delta\eta$, thus $u=\Delta^{-1}\delta\eta$. Therefore
\begin{equation}\label{eq:exactPe}
   \Pe\eta=du=d\Delta^{-1}\delta\eta.
\end{equation}
Any two potentials differ by a constant and give the same $du$.
The formula also shows that $\Pe$ preserves smoothness.

\begin{theorem}\label{lap:thm:form-estimate}
Assume \eqref{eq:transfer-gap}. Every smooth real differential
form $\chi=\sum_{p=0}^n\chi_p$ satisfies
\begin{equation}\label{lap:eq:form-estimate}
 \norm{\LA\chi}^2
 \ge\sum_{p=0}^n p(n-p)\norm{\chi_p}^2
     +\norm{\Pe\chi_1}^2.
\end{equation}
\end{theorem}

\begin{proof}
Under \eqref{lap:quadratic-gap}, the left-hand side of
\eqref{lap:eq:quadratic-identity} is nonnegative for every $u$.
Orthogonality of the exact projection $\Pe$ gives
\[
 \norm{du}^2+2\Ltwo{du}{\chi_1}
 =\norm{du+\Pe\chi_1}^2-\norm{\Pe\chi_1}^2.
\]
Choose $u=-\Delta^{-1}\delta\chi_1$, so that $du+\Pe\chi_1=0$. Substitution into \eqref{lap:eq:quadratic-identity}
proves \eqref{lap:eq:form-estimate}.
\end{proof}

\begin{lemma}\label[lemma]{lap:lem:kernel}
Under \eqref{eq:transfer-gap},
\begin{equation}\label{lap:eq:kernel}
 \ker\LA=\operatorname{span}_{\R}\{1,\volM\}.
\end{equation}
\end{lemma}

\begin{proof}
Let $\chi\in\ker\LA$. Elliptic regularity makes $\chi$ smooth,
and \eqref{lap:eq:form-estimate} forces $\chi_p=0$ for
$1\le p\le n-1$, an empty range when $n=1$.
Using \eqref{lap:eq:endpoint-A}, the kernel equation reduces to
\[
 0=\LA\chi=d\chi_0+\delta\chi_n.
\]
The two terms are exact and coexact, respectively, and hence
orthogonal even when $n=2$. Both must vanish. Connectedness then
implies that $\chi_0$ is constant and $\chi_n$ is a constant
multiple of the volume form $\volM$. Conversely, both $d+\delta$ and
$\KK$ annihilate $1$ and $\volM$, and so do $\LA$.
\end{proof}

\begin{proposition}\label[proposition]{lap:prop:inverse}
Under \eqref{eq:transfer-gap}, every smooth real mean-zero function
$f$ satisfies
\begin{equation}\label{lap:eq:inverse}
 \norm{f}^2\ge n\Ltwo{f}{\Delta^{-1}f}.
\end{equation}
\end{proposition}

\begin{proof}
Since $M$ is closed and $\LA$ is elliptic, its realization
$\LA:H^1\to L^2$ has closed range. By self-adjointness, the Fredholm
alternative \cite[Ch.~III, \S5, Theorem~5.5, p.~195]{LawsonMichelsohn}
therefore gives
\[
 \operatorname{ran}\LA
 =\bigl((\operatorname{ran}\LA)^\perp\bigr)^\perp
 =(\ker\LA^*)^\perp=(\ker\LA)^\perp.
\]
As a zero-form, $f$ satisfies $\Ltwo{f}{1}=0$ by the mean-zero
assumption and $\Ltwo{f}{\volM}=0$ by degree. Thus
\cref{lap:lem:kernel} gives $f\in(\ker\LA)^\perp=\operatorname{ran}\LA$,
so there is a form $\chi\in H^1$ satisfying
\begin{equation}\label{lap:eq:solve}
 \LA\chi=f.
\end{equation}
The coefficients of $\LA$ are real, so the solution can be
chosen real; elliptic regularity makes it smooth. The zero-form
component of \eqref{lap:eq:solve} is
\[
 \delta\chi_1=f.
\]
Therefore, by \eqref{eq:exactPe} the exact part of $\chi_1$ is 
\begin{equation}\label{lap:eq:exact-component}
 \Pe\chi_1=d\Delta^{-1}f.
\end{equation}
Retaining the degree-one term in \eqref{lap:eq:form-estimate}
and using $\norm{\Pe\chi_1}\le\norm{\chi_1}$, we obtain
\begin{align*}
 \norm{f}^2
 &=\norm{\LA\chi}^2\\
 &\ge(n-1)\norm{\chi_1}^2+\norm{\Pe\chi_1}^2\\
 &\ge n\norm{\Pe\chi_1}^2
 =n\norm{d\Delta^{-1}f}^2
 =n\Ltwo{f}{\Delta^{-1}f}.
\end{align*}
\end{proof}

\begin{proof}[Proof of \cref{thm:laplace-transfer}]
Let $f\ne0$ satisfy $\Delta f=\mu f$ with $\mu>0$. The function
$f$ has mean zero, so \cref{lap:prop:inverse} yields
\[
 \norm{f}^2\ge n\Ltwo{f}{\Delta^{-1}f}=\frac n\mu\norm{f}^2.
\]
Thus every positive eigenvalue satisfies $\mu\ge n$. The
coordinate identity \eqref{lap:coordinate} supplies the eigenvalue
$n$, proving $\lambda_1(\Delta)=n$.
\end{proof}

\medskip
\noindent\textbf{Acknowledgements and AI Assistance Disclosure.}
The first author J. Ge is grateful to Professor Xiao Zhang for his valuable discussions and enlightening online invited-short-course on Dirac operators in the year 2020. At that time, J. Ge proposed the assertions of the two results ($\lambda_1^+(D_M)=\frac n2$ for both normals $\pm\nu$ and this implies $\lambda_1(\Delta_M)=n$ in the isoparametric case) of the paper. Many thanks to ChatGPT for introducing the leafwise estimate and the dimension-unified doubled spinor bundle construction, which facilitate us in the proofs. 
The authors checked and wrote all math contents, and take full responsibility for this manuscript.

\begingroup
\linespread{1.0}\selectfont

\endgroup
\end{document}